\documentclass[12pt, reqno]{amsart}
\usepackage[a4paper, top=2.5cm, bottom=2.5cm, left=2cm, right=2cm]{geometry}

\usepackage[english]{babel}

\usepackage{amsmath, amsthm, amscd, amsfonts, amssymb, color}
\usepackage[bookmarksnumbered, colorlinks, plainpages]{hyperref}
\usepackage{float}
\usepackage{bbm}
\usepackage{mathrsfs}
\usepackage{dsfont}
\usepackage{csquotes}
\usepackage{latexsym}
\usepackage{exscale}

\usepackage{tikz}
\usepackage{tikz-cd}
\usetikzlibrary{shapes,arrows,positioning,calc}
\usetikzlibrary{matrix,arrows,decorations.pathmorphing}

\hypersetup{colorlinks=true,linkcolor=red, anchorcolor=green, citecolor=cyan, urlcolor=red, filecolor=magenta, pdftoolbar=true}

\newtheorem{theorem}{Theorem}[section]

\theoremstyle{definition}

\theoremstyle{remark}

\numberwithin{equation}{section}

\begin{document}

\setcounter{page}{1}

\title[Measuring the Infinite]{Measuring the Infinite: Interpolation Theory, Lorentz Spaces, and Dispersive PDEs}
 
\author[A. G. AKSOY, D. A. THIONG]{Asuman G\"{u}ven AKSOY$^{*}$, Daniel Akech Thiong$^{1}$}

\address{$^{*}$Department of Mathematics, Claremont McKenna College, 850 Columbia Avenue, Claremont, CA  91711, USA.}
\email{\textcolor[rgb]{0.00,0.00,0.84}{aaksoy@cmc.edu}}

\address{$^{1}$Department of Mathematics, Claremont Graduate University, 710 N. College Avenue, Claremont, CA  91711, USA.}
\email{\textcolor[rgb]{0.00,0.00,0.84}{daniel.akech@cgu.edu}}



\begin{abstract}
This article explores the vital role of interpolation theory and Lorentz spaces in the rigorous analysis of linear operators. While classical Lebesgue spaces ($L_{p}$) successfully measure the magnitude of functions, they frequently fail to bound evolution operators at critical endpoints of $p=1$ or $p = \infty$ because they conflate a function's amplitude with its spatial spread. To resolve this analytic bottleneck, we introduce distribution functions and decreasing rearrangements, culminating in the construction of Lorentz spaces ($L_{p, q}$). By utilizing the Complex (Riesz-Thorin), Real (Peetre's K-functional), and Marcinkiewicz methods of interpolation, these highly sensitive intermediate spaces act as geometric bridges between endpoint extremes. We conclude by applying this abstract framework to two distinct illustrative models: deriving the continuous smoothing decay of the parabolic Heat equation, and establishing the foundational dispersive Strichartz estimates for the dispersive free Schrödinger equation. 
\end{abstract} 

\maketitle

\section{Introduction}

The classical Lebesgue spaces, denoted as $L_p$, provide a versatile spectrum of mathematical measures for functional analysis. By integrating the $p$-th power of a function's absolute value, these spaces establish a complete normed vector geometry that is exceptionally well-behaved when strictly bounded between one and infinity. However, an inherent structural limitation of the $L_p$ norm is that it conflates two distinct geometric properties of a function: its absolute amplitude and its spatial spread. A function with a microscopic width but a massive amplitude can yield the exact same norm value as a function with a tiny amplitude stretching across the entire real line.

This rigidity creates severe analytic bottlenecks when dealing with operators that exhibit singular behavior at the critical endpoints of $p=1$ and $p=\infty$. To observe this mathematical failure directly, consider the classical Hardy operator, which serves as a fundamental averaging mechanism. For a non-negative integrable function $f$, the Hardy operator is defined as:
\[
    Hf(x) = \frac{1}{x} \int_0^x f(t) \, dt
\]
It is a well-established theorem that the Hardy operator maps $L_p$ into $L_p$ for $p > 1$. Specifically, the norm of the operator satisfies the strict inequality:
\[
    \|Hf\|_p \le \frac{p}{p-1} \|f\|_p
\]
The constant $\displaystyle \frac{p}{p-1}$ is the best possible bound for this operator (\cite{Hardy}, Theorem 327). As $p$ approaches $1$ from the right, the term $p-1$ approaches zero, and the bounding constant explodes to infinity. Therefore, $\|H\|_p \to \infty$, proving that the Hardy operator is unbounded on $L_1$. If an analyst only possessed the monolithic spaces $L_1$ and $L_\infty$ to work with, bounding singular integral operators would appear intractable.

The central goal of this exposition is to appeal to Interpolation Theory to address this limitation. Rather than treating interpolation as a mere topological curiosity, this article provides an expository journey through the construction of rearrangement-invariant structures. We will guide the reader through the decoupling of size and spread using distribution functions, the restoration of sub-additivity via Peetre's K-functional, and the foundational interpolation theorems of Riesz-Thorin and Marcinkiewicz. 

By replacing rigid Lebesgue spaces with finer, rearrangement-invariant structures, we gain the exact analytical control required to bound complex operators. Finally, as a capstone demonstration of this abstract machinery's utility, the paper will conclude by applying these techniques to two concrete problems in partial differential equations: the thermal diffusion of the heat equation and the derivation of the Strichartz estimates for the free Schrödinger equation.

\section{Rearrangement-Invariant Spaces}

\subsection{The Distribution Function and Decreasing Rearrangements}

To systematically decouple the size and spread of a measurable function $f$, we ignore the specific locations where the function achieves its values and instead measure the size of its level sets. Let $(X, \mu)$ be a $\sigma$-finite measure space. For each measurable function $f$, we define the \textit{distribution function} $d_f$ by:
\[
    d_f(\alpha) = \mu(\{x \in X : |f(x)| > \alpha\})
\]
The distribution function returns the total measure (or ``width'') of the domain where the function's absolute value exceeds a given height $\alpha$. It is a non-increasing, right-continuous function. By utilizing the distribution function, we strip away the oscillating, potentially chaotic spatial layout of $f$ and retain only its volumetric profile.

However, $d_f$ acts on the range of $f$. To create a standardized function defined on the domain $\mathbb{R}^+$, we construct the \textit{decreasing rearrangement} $f^*$, defined as:
\[
    f^*(t) = \inf \{\alpha > 0 : d_f(\alpha) \le t\}
\]
It conceptually takes all the erratic, jagged values of $f$ across the domain, sorts them from largest to smallest, and pushes them tightly against the y-axis, creating a right-continuous, non-increasing profile. A crucial property of this transformation is that functions $f$ and $f^*$ are equimeasurable; their distribution functions coincide. Consequently, if $f \in L_p$, then $f^* \in L_p$, and their norms are identical: $\|f\|_p = \|f^*\|_p$ (\cite{Berg}, \cite{DeVore}, \cite{Grafakos}, \cite{SteinWeiss}).

\subsection{The Role of the Hardy Operator}

While the decreasing rearrangement provides beautiful geometric clarity, it harbors a complication for functional analysis: it is not sub-additive. In general, the rearrangement of a sum is not bounded by the sum of the rearrangements; that is, $(f+g)^* \le f^* + g^*$ is generally false. 

To conduct rigorous functional analysis and solve differential equations, we require a complete normed space, such as a Banach space. Without the triangle inequality, we merely have a quasi-norm, which makes proving operator bounds exceptionally difficult. We must artificially induce sub-additivity into the rearrangement profile.

We rectify this by applying the aforementioned Hardy operator (introduced in Section 1) to compute the running integral average of the decreasing rearrangement. We denote this maximal function as $f^{**}$:
\[
    f^{**}(t) = \frac{1}{t} \int_0^t f^*(s) \, ds
\]
Because $f^{**}$ represents an integral average of a non-increasing function $f^*$, it is inherently larger than or equal to $f^*$. More importantly, the integration process restores sub-additivity. The functional strictly satisfies $(f+g)^{**} \le f^{**} + g^{**}$ (\cite{Berg}, \cite{Ben}, \cite{DeVore}, \cite{Hardy}). With this sub-additive, rearrangement-invariant profile in hand, we are equipped to define a highly sensitive class of intermediate function spaces.

\subsection{Lorentz Spaces}

Using the rearrangement profiles, we define the \textit{Lorentz spaces}, denoted $L_{p,q}$. These spaces act as a vast, two-dimensional generalization of the one-dimensional Lebesgue scale (\cite{Berg}, \cite{DeVore}). 

Let $(X, \mu)$ be a measure space. For $1 \le p < \infty$ and $1 \le q < \infty$, the Lorentz space $L_{p,q}$ is canonically defined as the set of all measurable functions $f$ such that:
\[
    \|f\|_{L_{p,q}} = \left( \frac{q}{p} \int_0^\infty \left[t^{1/p} f^*(t)\right]^q \, \frac{dt}{t} \right)^{1/q} < \infty
\]
For the endpoint case where $q = \infty$, the space $L_{p,\infty}$ is defined by the supremum:
\[
    \|f\|_{L_{p,\infty}} = \sup_{t>0} t^{1/p} f^*(t) < \infty
\]

While the functional defined above using the decreasing rearrangement $f^*$ recovers the Lebesgue space $L_p$ precisely when $p=q$, it merely constitutes a quasi-norm since $f^*$ is not sub-additive. To establish a true Banach space geometry equipped with a rigorous norm, functional analysts routinely replace $f^*$ with the averaged maximal function $f^{**}$. For $1 < p \le \infty$, the Hardy operator is bounded, and the resulting $f^{**}$ functional yields an equivalent, sub-additive norm.

However, as noted in Section 1, the Hardy operator is unbounded at $p=1$. Consequently, integrating $f^{**}$ for $p=1$ introduces a non-integrable singularity. Therefore, to correctly recover $L_1$ and maintain consistent behavior at this critical endpoint, one must rigorously decouple the cases and rely on the quasi-norm $f^*$ rather than $f^{**}$ precisely when $p=1$ \cite{SteinWeiss}.

The parameters $p$ and $q$ serve distinct, highly specialized roles. The primary exponent $p$ dictates the overarching global size and scaling of the function, identically to its role in a standard Lebesgue space. The secondary exponent $q$ provides fine-tuned control over the local behavior and the asymptotic decay tail of the function's rearrangement. For a fixed primary exponent $p$, altering $q$ creates a nested hierarchy of spaces. It can be shown that if $q_1 < q_2$, then $L_{p,q_1} \subset L_{p,q_2}$.

The space $L_{p,1}$ is the smallest and most restrictive Lorentz space for a given $p$, demanding that the function's peaks and tails be highly controlled. Conversely, the space $L_{p,\infty}$, often called the ``weak space,'' is the largest and most forgiving, capturing functions that just barely fail to be integrable in standard $L_p$ due to logarithmic divergences (\cite{Ben},\cite{Berg}, \cite{SteinWeiss},\cite{Grafakos}).

\section{The Bridges of Functional Analysis: Interpolation Theorems}

Proving that a complex differential operator maps one highly specific space into another (e.g., $L_p \to L_q$) through brute-force integration is often mathematically intractable. Interpolation theory provides an elegant structural shortcut. If we can prove that an operator behaves well at two relatively simple ``extreme'' endpoint cases, interpolation theorems act as geometric bridges, automatically guaranteeing that the operator is bounded on all the intermediate spaces existing between those endpoints.

To formalize this, we define the concept of a \textit{Banach couple}. A pair of Banach spaces $(X_0, X_1)$ is a Banach couple if both $X_0$ and $X_1$ are continuously embedded into a common, larger Hausdorff topological vector space $V$. An intermediate space $X$ is any normed space residing securely between these boundaries: $X_0 \cap X_1 \subset X \subset X_0 + X_1$.

An intermediate space $X$ is classified specifically as an \textit{interpolation space} if any linear operator that maps $X_0$ boundedly to itself, and $X_1$ boundedly to itself, is automatically guaranteed to map $X$ boundedly to itself. 

\subsection{The Complex Method and The Riesz-Thorin Theorem}

The complex method of interpolation is deeply intertwined with the theory of analytic functions and the maximum modulus principle \cite{Grafakos}. The bound $M_\theta \le M_0^{1-\theta} M_1^\theta$ we see in the following Riesz-Thorin Theorem is produced directly by the Hadamard three-lines theorem from complex analysis. 

\begin{theorem}[The Riesz-Thorin Interpolation Theorem]
Let $(X, \mu)$ and $(Y, \nu)$ be $\sigma$-finite measure spaces. Assume that $1 \le p_0, p_1, q_0, q_1 \le \infty$, and let $T$ be a linear operator defined on the simple functions. Suppose that $T$ is bounded at the two extreme endpoints:
\[
    T: L_{p_0} \to L_{q_0} \quad \text{with operator norm } M_0
\]
\[
    T: L_{p_1} \to L_{q_1} \quad \text{with operator norm } M_1
\]
Then, for any parameter $\theta \in (0,1)$, the operator $T$ uniquely extends to a bounded linear map on the intermediate space:
\[
    T: L_{p_\theta} \to L_{q_\theta}
\]
Furthermore, the operator norm $M_\theta$ in this intermediate space is  bounded by the geometric mean of the endpoint norms:
\[
    M_\theta \le M_0^{1-\theta} M_1^\theta
\]
provided that the intermediate exponents $p_\theta$ and $q_\theta$ are defined by the convex combinations:
\[
    \frac{1}{p_\theta} = \frac{1-\theta}{p_0} + \frac{\theta}{p_1}, \quad \frac{1}{q_\theta} = \frac{1-\theta}{q_0} + \frac{\theta}{q_1}
\]
\end{theorem}

The Riesz-Thorin theorem is a foundational result of log-convexity \cite{Berg}. However, its reliance on complex analyticity renders it structurally rigid; it applies exclusively to strictly linear operators and standard $L_p$ spaces \cite{Grafakos}. To handle sub-linear operators and to incorporate the nuanced tail controls of Lorentz spaces, we must turn to the real and weak methods.

\subsection{The Real Method: Peetre's K-Functional}

While the complex method relies on holomorphic functions, the real method relies on the pure real-variable geometry of decreasing rearrangements and magnitude optimization. Formulated primarily by Jaak Peetre and Jacques-Louis Lions, the real method relies on measuring the absolute optimal way to decompose a function into two parts.

For a function $f$ residing in the sum of a Banach couple $X_0 + X_1$, and for any time-like parameter $t > 0$, Peetre defined the \textit{K-functional} as:
\[
    K(f, t) = \inf_{f = f_0 + f_1} (\|f_0\|_{X_0} + t\|f_1\|_{X_1})
\]
The K-functional acts as a mathematical penalty function.  It measures the minimal cost required to split $f$ into a ``good'' bounded part $f_0 \in X_0$ and an unbounded, residual part $f_1 \in X_1$, where the ``bad'' part is penalized by the weight $t$. 

By integrating this K-functional against a weighted measure, we generate an entirely new parameterized family of real interpolation spaces, denoted $(X_0, X_1)_{\theta, q}$. The exact norm is defined by:
\[
    \|f\|_{\theta,q} = \left( \int_0^\infty (t^{-\theta} K(f,t))^q \, \frac{dt}{t} \right)^{1/q}
\]
A significant advantage of the real interpolation method is its deep, fundamental connection to Lorentz spaces. If we choose our starting Banach couple to be the ultimate endpoints of classical analysis---the space of integrable functions $L_1$ and the space of bounded functions $L_\infty$---and compute the real interpolation spaces between them, an astonishing equivalence emerges.

It can be rigorously proven that for $f \in L_1 + L_\infty$, the K-functional for the couple $(L_1, L_\infty)$ is exactly equal to the integrated decreasing rearrangement: $K(f, t) = \int_0^t f^*(s) \, ds$. Substituting this equivalence into the K-functional integral norm perfectly reproduces the definition of the Lorentz space $L_{p,q}$. More precisely, $$(L_1, L_\infty)_{\theta, q}= L_{p,q},\,\,\,  \mbox{where} \,\,\, \frac{1}{p}=1-\theta, \,\,\, 0< \theta<1, \,\, 1\leq q\leq \infty,$$ with equivalent norms, see \cite{Berg}, \cite{Triebel}, \cite{Lunardi}.

\subsection{The Marcinkiewicz Interpolation Theorem}
While the Riesz-Thorin theorem is elegant, it strictly requires strong operator bounds at the endpoints. The Marcinkiewicz Interpolation Theorem  \cite{Ben} brilliantly circumvents this restriction by utilizing the weak Lorentz spaces $L_{p,\infty}$, establishing the final pillar of classical interpolation theory.

\begin{theorem}[Marcinkiewicz Interpolation Theorem]
Let $1 \le p_0 < p_1 \le \infty$ and $1 \le q_0, q_1 \le \infty$ with $q_0 \neq q_1$. Suppose $T$ is a sub-linear operator bounded from $L_{p_0}$ to the weak space $L_{q_0,\infty}$ and from $L_{p_1}$ to the weak space $L_{q_1,\infty}$. Then for any $\theta \in (0,1)$, $T$ extends to a bounded operator between the strong Lebesgue spaces:
\[
    T: L_{p_\theta} \to L_{q_\theta}
\]
where $\frac{1}{p_\theta} = \frac{1-\theta}{p_0} + \frac{\theta}{p_1}$ and $\frac{1}{q_\theta} = \frac{1-\theta}{q_0} + \frac{\theta}{q_1}$.
\end{theorem}

This theorem solidifies the absolute necessity of Lorentz spaces: by allowing functions to reside in the slightly larger, more forgiving $L_{p,\infty}$ space at the extremes, we can still deduce perfect, strong $L_p$ bounds in the intermediate spaces. Interpolation theory grants us a complete map of the function space universe, and we are now ready to apply this map.

\section{Application to Thermal Diffusion and the Heat Equation}

We now pivot to applied physical mechanics as a capstone demonstration of our abstract frameworks. Consider a physical medium where an initial distribution of temperature, denoted $u(x,0) = f(x)$, is introduced at time zero. Our goal is to rigorously prove that the evolution operator governing heat diffusion maps rough, chaotic initial data into perfectly smooth, restricted spaces.

\subsection{The Heat Semigroup and Kernel Convolution}

The initial value Cauchy problem is defined as:
\[
    \begin{cases} 
        \partial_t u = \Delta u \\ 
        u(x,0) = f(x) 
    \end{cases}
\]
The solution at any time $t$ can be elegantly expressed as a spatial convolution of the initial data with the Gaussian heat kernel $K_t(x)$:
\[
    K_t(x) = \frac{1}{(4\pi t)^{n/2}} e^{-\frac{|x|^2}{4t}}
\]
Evaluating this Gaussian integral over $L_r(\mathbb{R}^n)$ yields a precise metric of how the kernel flattens over time:
\[
    \|K_t\|_r = C_r t^{-\frac{n}{2}\left(1 - \frac{1}{r}\right)}
\]

\subsection{Young's Inequality and Lorentz Space Refinements}

By substituting our precise evaluation of the heat kernel's norm into Young's inequality, we obtain the fundamental endpoint estimates for the heat semigroup $T_t$:
\[
    \|T_t f\|_p \le C t^{-\frac{n}{2}\left(\frac{1}{q} - \frac{1}{p}\right)} \|f\|_q
\]

By invoking the Riesz-Thorin interpolation theorem on the operator $T_t$, we readily obtain the boundedness criteria for the entire infinite spectrum of intermediate Lebesgue spaces $L_{q_\theta}$ and $L_{p_\theta}$ simultaneously:
\[
    \|T_t f\|_{p_\theta} \le C_\theta t^{-\alpha_\theta} \|f\|_{q_\theta}
\]

Standard Lebesgue spaces fail to capture the behavior of solutions when the initial data $f(x)$ possesses severe polynomial singularities, yet they reside perfectly within the weak Lorentz space $L_{p,\infty}$. By applying the real interpolation method (Peetre's K-functional) to the heat operator endpoints, mathematicians prove that $T_t$ extends continuously to a bounded operator between Lorentz spaces:
\[
    T_t : L_{q, \infty} \to L_{p, \infty}
\]
Through real interpolation, we rigorously confirm that the smoothing mechanism of the heat operator is so robust that it diffuses even borderline singular data (weak spaces) into highly regularized thermal states.

\section{Application to Quantum Dispersion and the Schrödinger Equation}

Unlike thermal systems, quantum wave packets do not merely ``smooth out''; their evolution is characterized by high-frequency wave components traveling at strictly faster phase velocities than low-frequency components \cite{TaoBook}. Over time, this differential velocity causes the localized wave packet to unravel across the spatial domain---a phenomenon known as dispersion \cite{TaoBook}. 

\subsection{The Unitary Group and the Absence of Damping}

Consider the linear initial value problem for the free Schrödinger equation in $\mathbb{R}^n$:
\[
    \begin{cases} 
        i\partial_t u + \Delta u = 0 \\ 
        u(x,0) = f(x) 
    \end{cases}
\]
Applying the Fourier Transform reveals the explicit integral formula for the Schrödinger evolution operator $T_t = e^{it\Delta}$:
\[
    u(x,t) = \frac{1}{(4\pi i t)^{n/2}} \int_{\mathbb{R}^n} e^{i\frac{|x-y|^2}{4t}} f(y) \, dy
\]
Because the evolution operator forms a unitary group on the Hilbert space $L_2$, the $L_2$ norm of the wave function is perfectly conserved:
\[
    \|u(\cdot, t)\|_2 = \|f\|_2
\]
This is our first endpoint estimate: perfect conservation. To mathematically capture the fading amplitude of dispersion, we look to the space of maximum heights, $L_\infty$. This establishes our second endpoint estimate---a strict decay bound:
\[
    \|u(\cdot, t)\|_\infty \le \frac{1}{(4\pi |t|)^{n/2}} \|f\|_1
\]

\subsection{The Dispersive Strichartz Estimates via Interpolation}

By applying the Riesz-Thorin complex interpolation theorem between the bounded norm at $p=2$ and the bounded norm at $p=\infty$, we yield the seminal \textit{dispersive estimate}:
\[
    \|e^{it\Delta} f\|_p \le C |t|^{-n\left(\frac{1}{2} - \frac{1}{p}\right)} \|f\|_{p'}
\]

\textbf{The Culmination in Lorentz Spaces:} Restricting the Schrödinger analysis to Riesz-Thorin interpolation only reveals half the picture. A central pillar of modern dispersive PDE theory relies on \textit{Strichartz estimates}, which represent mixed space-time bounds that are strictly required to prove the stability of complex non-linear Schrödinger variants \cite{Tao}. 

These profound Strichartz estimates are derived explicitly by utilizing the real interpolation methodology and the Lorentz space structures developed in Section 2. A pair of exponents $(q, r)$ is called \textit{admissible} if:
\[
    \frac{2}{q} + \frac{n}{r} = \frac{n}{2}
\]
where $2 \le q, r \le \infty$ and $(q, r, n) \neq (2, \infty, 2)$. The foundational Strichartz estimate guarantees the existence of a constant $C$ such that:
\[
    \|e^{it\Delta} f\|_{L_q(\mathbb{R}; L_r(\mathbb{R}^n))} \le C \|f\|_{L_2(\mathbb{R}^n)}
\]

By interpolating the Schrödinger operator across time using Peetre's abstract K-functional, mathematicians generate these mixed bounds mapping initial data into the Bochner space $L_q(\mathbb{R}; L_r(\mathbb{R}^n))$, a natural extension of Lebesgue integral to functions taking values in a Banach space, which measures how the spatial $L_r$ norm of the wave packet varies over a timeline measured by $L_q$). Through the highly tuned, tail-controlling geometry of Lorentz spaces, the sharpest bounds required for quantum analysis become mathematically transparent. The original restriction theorem appears in \cite{Str}; see also \cite{Tao} for endpoint estimates and \cite{Cazenave,Evans} for further PDE applications.

\begin{table}[htbp]
\centering
\caption{Comparison of Boundedness Criteria for PDE Operators derived via Interpolation Theory.}
\renewcommand{\arraystretch}{1.5}
\begin{tabular}{p{0.15\textwidth} p{0.2\textwidth} p{0.2\textwidth} p{0.15\textwidth} p{0.2\textwidth}}
\hline
\textbf{PDE System} & \textbf{Physical Mechanism} & \textbf{Evolution Operator} & \textbf{Endpoint Limits} & \textbf{Interpolated Estimate Bound} \\
\hline
Heat Equation & Parabolic / Dissipative & Convolution Semigroup & $L_1 \to L_1$,\newline $L_1 \to L_\infty$ & $\|T_t f\|_{p} \lesssim t^{-\frac{n}{2}(\frac{1}{q} - \frac{1}{p})} \|f\|_{q}$ \\
Schrödinger Eq. & Dispersive & Unitary Group & $L_2 \to L_2$,\newline $L_1 \to L_\infty$ & $\|e^{it\Delta} f\|_{p} \lesssim t^{-n(\frac{1}{2} - \frac{1}{p})} \|f\|_{p'}$ \\
\hline
\end{tabular}
\end{table}

\section{Visualizing the Abstraction}

To solidify the abstract mechanics of interpolation and function spaces discussed in the preceding sections, we present three visual representations. A well-constructed visual metaphor reduces cognitive load and illuminates the underlying geometry that abstract algebra alone often obscures. 

\subsection{The Decreasing Rearrangement Profile}

Figure 1 illustrates the foundational mechanics of the decreasing rearrangement $f^*$, as defined in Section 3. The left panel displays a highly oscillatory, non-negative function $f(x) = e^{-x^{2}}\sin ^{2} (10x) + \frac{1}{2} e^{-(x-2)^{2}}$ across the real number line, exhibiting jagged, asymmetric peaks. For an arbitrary height $\alpha$, the shaded segments represent the distribution measure $d_f(\alpha)$---the total width where the curve exceeds $\alpha$. 

The right panel displays the function $f^*$. The rearrangement operator acts as a mathematical sorting algorithm: it preserves the exact area of the level sets (and therefore all $L_p$ integrals) but completely destroys the original spatial coordinates, sweeping all the function's mass from highest density to lowest density to form a monotonically decreasing curve.

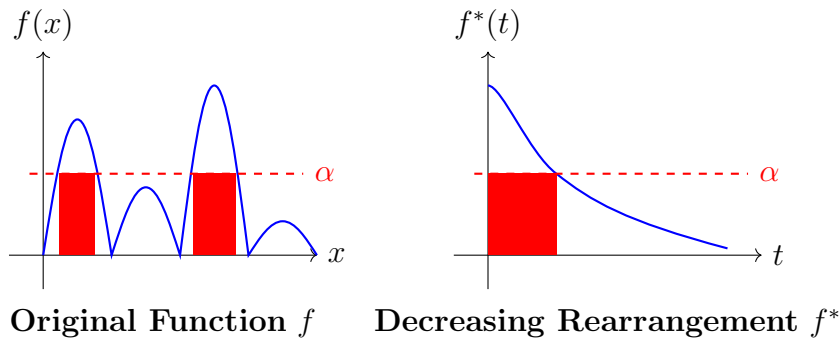
\begin{figure}[htbp]
\centering
\begin{tikzpicture}[scale=0.9]
    \draw[->] (-0.5,0) -- (4,0) node[right] {$x$};
    \draw[->] (0,-0.5) -- (0,3) node[above] {$f(x)$};
    \draw[thick, blue] (0,0) sin (0.5,2) cos (1,0) sin (1.5,1) cos (2,0) sin (2.5,2.5) cos (3,0) sin (3.5, 0.5) cos (4,0);
    
    \draw[dashed, red, thick] (-0.2, 1.2) -- (3.8, 1.2) node[right] {$\alpha$};
    
    \fill[red, opacity=0.3] (0.24,0) rectangle (0.76, 1.2); 
    \fill[red, opacity=0.3] (2.19,0) rectangle (2.81, 1.2);
    
    \node at (1.75,-1) {\textbf{Original Function} $f$};

    \begin{scope}[xshift=6.5cm]
    \draw[->] (-0.5,0) -- (4,0) node[right] {$t$};
    \draw[->] (0,-0.5) -- (0,3) node[above] {$f^*(t)$};
    
    \draw[thick, blue] (0,2.5) .. controls (0.2,2.5) and (0.6,1.5) .. (1.0, 1.2) .. controls (1.5, 0.8) and (2,0.5) .. (3.5,0.1);
    
    \draw[dashed, red, thick] (-0.2, 1.2) -- (3.8, 1.2) node[right] {$\alpha$};
    
    \fill[red, opacity=0.3] (0,0) rectangle (1.0, 1.2);
    
    \node at (1.75,-1) {\textbf{Decreasing Rearrangement} $f^*$};
    \end{scope}
\end{tikzpicture}
\caption{The Decreasing Rearrangement. The total measure of the domain where the function exceeds height $\alpha$ (red shaded regions) is preserved and consolidated against the y-axis, decoupling amplitude from spatial spread.}
\end{figure}

\subsection{The Riesz Convexity Square}

The Riesz-Thorin interpolation theorem is classically illustrated using the convexity square (Figure 2). The diagram plots inverse exponents on a $[0,1] \times [0,1]$ coordinate system. The x-axis represents the input space index $1/p$, and the y-axis represents the output space index $1/q$. The origin corresponds to the extreme boundary $L_\infty \to L_\infty$, while $(1,1)$ corresponds to $L_1 \to L_1$.

The points $A$ and $B$ represent two known, explicitly proven endpoint bounds of a given operator $T$. The Riesz convexity square translates complex algebraic fractions into a geometric truth: if an operator is bounded at two isolated locations within the function space universe, it is automatically bounded everywhere along the linear trajectory connecting them (Point $C_\theta$).

\begin{figure}[htbp]
    \centering
    \begin{tikzpicture}[scale=4]
        \draw[step=0.5, gray!40, very thin, dashed] (0,0) grid (1,1);
        \draw[thick] (0,0) rectangle (1,1);

        \node[below] at (0.5, -0.05) {Input Space Index ($1/p$)};
        \node[left] at (-0.05, 0.5) {Output Space Index ($1/q$)};

        \coordinate (A) at (0.2, 0.8);
        \coordinate (B) at (0.9, 0.3);

        \draw[very thick, blue] (A) -- (B);

        \coordinate (C) at ($(A)!0.45!(B)$);

        \filldraw[black] (A) circle (0.6pt) node[above right] {$A(1/p_0, 1/q_0)$};
        \filldraw[black] (B) circle (0.6pt) node[above right] {$B(1/p_1, 1/q_1)$};
        \filldraw[red] (C) circle (0.8pt) node[above right, text=red] {$C_\theta(1/p_\theta, 1/q_\theta)$};

        \node[below left] at (0,0) {$0$};
        \node[below] at (1,0) {$1$};
        \node[left] at (0,1) {$1$};
    \end{tikzpicture}
    \caption{The Riesz Convexity Square. Interpolation provides intermediate operator bounds precisely along the convex combination of known endpoints.}
\end{figure}
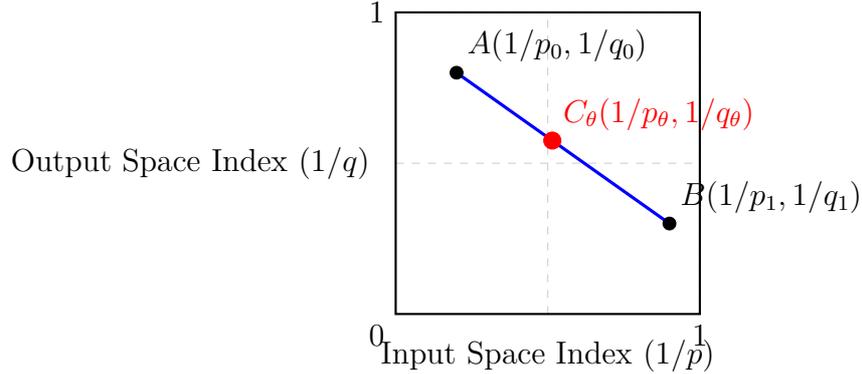

\subsection{Diffusive vs. Dispersive Evolution}

Figure 3 physically differentiates the abstract operator bounds derived for the Heat and Schrödinger equations. The top sequence models the Heat equation: an initial localized spike simply flattens and widens symmetrically. The amplitude decays uniformly, and the gradient of the peak vanishes into a smooth bell curve, representing pure thermal dissipation.

Conversely, the bottom sequence models the Schrödinger equation. While the overall bounding envelope (dashed red) widens and decays in maximum amplitude $L_\infty$, the internal structure exhibits chaotic, high-frequency oscillations. This visualizes \textit{why} the equations behave differently under interpolation: heat provides a highly regularized map into strictly smooth spaces, whereas quantum evolution merely decays the envelope while preserving the integral square energy of the internal ripples.

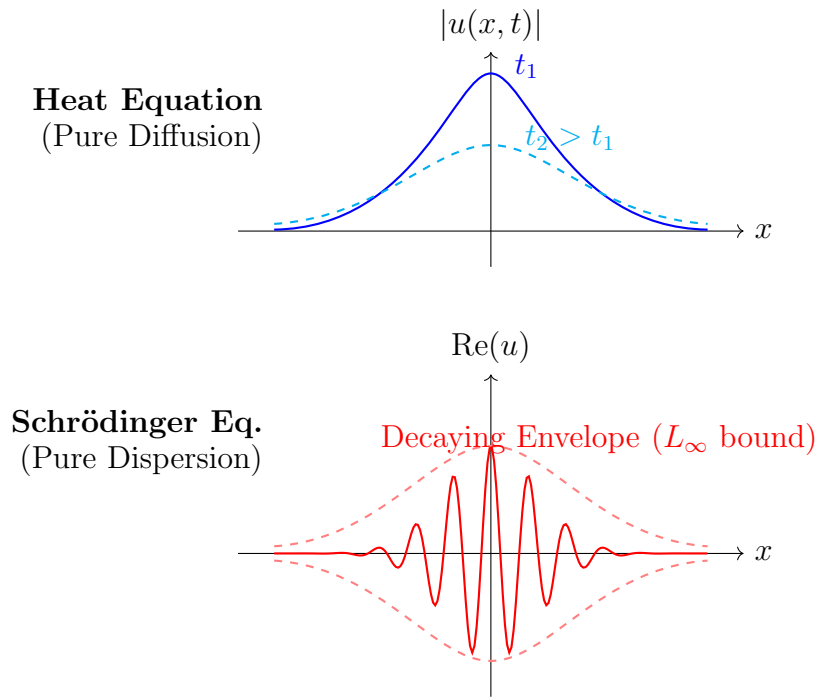
\begin{figure}[htbp]
\centering
\begin{tikzpicture}[scale=0.95]
    \begin{scope}[yshift=4.5cm]
    \node[left] at (-3, 1.8) {\textbf{Heat Equation}};
    \node[left] at (-3, 1.3) {(Pure Diffusion)};
    
    \draw[->] (-3.5,0) -- (3.5,0) node[right] {$x$};
    \draw[->] (0,-0.5) -- (0,2.5) node[above] {$|u(x,t)|$};
    
    \draw[thick, blue] (-3,0.02) .. controls (-1,0.05) and (-0.5,2.2) .. (0,2.2) .. controls (0.5,2.2) and (1,0.05) .. (3,0.02);
    \node[blue] at (0.5, 2.3) {$t_1$};
    
    \draw[thick, cyan, dashed] (-3,0.1) .. controls (-1.5,0.2) and (-0.8,1.2) .. (0,1.2) .. controls (0.8,1.2) and (1.5,0.2) .. (3,0.1);
    \node[cyan] at (1.1, 1.3) {$t_2 > t_1$};
    \end{scope}

    \begin{scope}[yshift=0cm]
    \node[left] at (-3, 1.8) {\textbf{Schrödinger Eq.}};
    \node[left] at (-3, 1.3) {(Pure Dispersion)};
    
    \draw[->] (-3.5,0) -- (3.5,0) node[right] {$x$};
    \draw[->] (0,-2.0) -- (0,2.5) node[above] {$\text{Re}(u)$};
    
    \draw[thick, red!50, dashed] (-3,0.1) .. controls (-1.5,0.2) and (-0.8,1.5) .. (0,1.5) .. controls (0.8,1.5) and (1.5,0.2) .. (3,0.1);
    \draw[thick, red!50, dashed] (-3,-0.1) .. controls (-1.5,-0.2) and (-0.8,-1.5) .. (0,-1.5) .. controls (0.8,-1.5) and (1.5,-0.2) .. (3,-0.1);
    
    \draw[thick, red, samples=200, domain=-3:3] plot (\x, {1.5 * exp(-1.2*\x*\x) * cos(deg(12*\x))});
    
    \node[red] at (1.5, 1.6) {Decaying Envelope ($L_\infty$ bound)};
    \end{scope}
\end{tikzpicture}
\caption{Diffusive vs. Dispersive Evolution. Top: Thermal diffusion strictly flattens the amplitude over time. Bottom: Quantum dispersion decays the bounding envelope while high-frequency internal oscillations spread the probability wave packet.}
\end{figure}

\newpage
\section{Conclusion}
The intersection of abstract interpolation theory and the partial differential equations of mathematical physics yields a profoundly elegant narrative. Historically, the meticulous construction of rearrangement-invariant function spaces and the K-functional were viewed as exercises in pure topological generalization. However, as this exposition demonstrates, these tools form the exact analytical infrastructure required to quantify bounding behavior that Lebesgue spaces cannot resolve. By shifting our focus from pure integrations to the highly sensitive geometric controls of Lorentz spaces, Riesz-Thorin, and Marcinkiewicz interpolation, we gain the definitive mathematical language required to measure the infinite.

\bibliographystyle{amsplain}

\end{document}